\documentclass[11pt]{article}

\usepackage{amsfonts,amsmath,latexsym}
\usepackage{amstext,amssymb}
\usepackage{amsthm}
\usepackage{color}
\usepackage{soul}

\usepackage{tikz}

\newcommand {\R}{\mathbb{R}}

\newcommand {\Sect}{\operatorname{Sect}}

\newtheorem {thm} {Theorem}

\newtheorem {lemma}[thm]{Lemma}
\newtheorem {cor} [thm] {Corollary}

\newcommand{\beq}{\begin{equation}}
\newcommand{\eeq}{\end{equation}}
\numberwithin{equation}{section}

\begin {document}

\title {Monotonicty of the first Dirichlet eigenvalue of the Laplacian on manifolds of non-positive curvature} 
\author{Tom Carroll \footnote {School of Mathematical 
Sciences, University College Cork, Ireland, {\tt t.carroll@ucc.ie}}
and Jesse Ratzkin \footnote{Department of Mathematics and Applied 
Mathematics, University of Cape Town, South Africa 
{\tt jesse.ratzkin@uct.ac.za}}}

\maketitle 

\begin {abstract} Let $(M,g)$ be a complete Riemannian manifold with nonpositive 
scalar curvature, let $\Omega \subset M$ be a suitable domain, and 
let $\lambda(\Omega)$ be the first Dirichlet eigenvalue of the Laplace-Beltrami 
operator on $\Omega$. We prove several bounds for the rate of decrease of 
$\lambda(\Omega)$ as $\Omega$ increases, and a result comparing the 
rate of decrease of $\lambda$ before and after a conformal diffeomorphism. 
Along the way, we prove a 
reverse-H\"older inequality for the first eigenfunction, which generalizes results 
of Chiti to the manifold setting and maybe be of independent interest. \end {abstract} 

\section {Introduction} 

Let $(M,g)$ be a complete Riemannian manifold, and let $\Omega
\subset M$ be a domain with $\bar \Omega$ compact and $\partial 
\Omega \in \mathcal{C}^\infty$. The first Dirichlet eigenvalue of the 
Laplace-Beltrami operator on $\Omega$ is a natural and important 
object. It controls the 
slowest rate of heat dissipation from $\Omega$, the largest value 
of the expected exit time of a Brownian particle from $\Omega$, and 
the fundamental frequency of vibration of $\Omega$ (when considered 
as a vibrating membrane with stationary boundary). Many results have
linked the first eigenvalue and its associated 
eigenfunction to the geometry of $\Omega$, and also to that of the 
ambient space $(M,g)$. 

The first Dirichlet eigenvalue $\lambda (\Omega)$ of the 
Laplace-Beltrami $\Delta_g$ operator on $\Omega$ is 
$$
\lambda(\Omega) = \inf \left \{ \frac{\int_\Omega |\nabla u|^2 dm}
{\int_\Omega u^2 \,dm} : u \in W^{1,2}_0(\Omega) \right \}, 
$$ 
where $dm$ represents the volume element on $M$. 
This infimum is realized by a nontrivial function $\phi$, which satisfies 
$$\Delta_g \phi + \lambda \phi = 0, \quad \left. \phi \right |_{\partial 
\Omega} = 0, \quad \phi> 0 \textrm{ inside }\Omega.$$
If $\bar \Omega$ is compact then $\lambda$ is a positive, simple eigenvalue. 
In the case that $(M,g)$ is Euclidean space, a fundamental 
result is the Faber-Krahn inequality: 
$$\lambda(\Omega) \geq \lambda (B_1) \left ( \frac{\omega_n}
{|\Omega|_0} \right )^{2/n},$$
and equality can only occur if $\Omega$ is a round ball. Here $\omega_n$ 
is the volume of the unit ball in $\R^n$ and $|\Omega|_0$ is the volume 
of $\Omega$. Similar inequalities hold if the ambient space is hyperbolic
space, or the ambient space is a sphere and $\Omega$ is convex. One can find an excellent 
introduction to all this and more in Chapter 1 of \cite{Chav2} 

Our main goal is to explore the rate of change of 
$\lambda(\Omega_t)$, where $\Omega_t$ is an evolving family 
of domains. More precisely, we let $\zeta(t,p)$ be a one-parameter 
family of diffeomorphisms, and let $\Omega_t = \zeta(t,\Omega)$. 
Suppose that
$\left \langle \eta, \frac{\partial \zeta}{\partial t} \right \rangle > 0,$
where $\eta$ is the outward unit normal of $\Omega_t$. 
Then $\Omega_s \subset \Omega_t$ for $s<t$, so 
domain monotonicty implies $\lambda$ is a decreasing function of 
$t$, and we would like to estimate its rate of decrease. 

The isoperimetric inequality is a key tool we use, so we will 
need some standing
hypotheses on the ambient space $(M,g)$ and the domain $\Omega$. 
We always assume one of the following holds:
\begin {itemize} 
\item either $(M,g)$ is compact, and $S_g \leq -n(n-1)\kappa^2$ for 
some fixed $\kappa \in \R$, and $|\Omega|_g$ is sufficiently small,
\item or $(M,g)$ is complete, and $S_g \leq -n(n-1)\kappa^2$ for 
some fixed $\kappa \in \R$, and $\Omega$ is contained in a 
small geodesic ball $B_r$, whose radius might depend on its 
center,
\item or $dim(M) = 2$ and $S_g \leq 0$ is nonpositive. 
\end {itemize} 
Here $S_g$ is the scalar curvature of the metric $g$. 
One should be able to prove  results 
similar to ours in the presence of a weaker scalar curvature 
bound, such as $S_g \leq n(n-1)\kappa^2$, using straight-forward
adaptations of our proofs below. 

Three examples of what we are able to prove are the 
following theorems. 
\begin {thm} \label{evolution-thm1} 
Let $(M,g)$ be a complete Riemannian manifold, and let $\Omega 
\subset M$ be a domain with $\bar \Omega$ compact and $\partial 
\Omega \in \mathcal{C}^\infty$. Suppose that $\partial \Omega$ 
evolves with velocity vector $\eta$, the unit length outward normal to 
$\Omega$, and let $\lambda(t) = \lambda (\Omega_t)$. 
If $\dim (M) = 2$ and $(M,g)$ has nonpositive Gauss curvature then 
$$\frac{d}{dt} \log (\lambda) \leq -\frac{4\pi}{|\partial \Omega|_g}.$$
If $\dim (M) =n\geq 3$, and $(M,g)$ has nonpositive scalar curvature and 
$\Omega$ is sufficiently small (see Section \ref{isop-sec} below)
then 
$$\frac{d}{dt} \left ( \lambda^{\frac{n-2}{2} }\right ) \leq -\left ( \frac{n-2}{2} 
\right ) \frac{K}{|\partial \Omega|_g},$$
where $K$ is a constant depending only on $n$. Equality in either case 
can only occur if $\Omega$ is isometric to a round ball in the appropriate 
dimensional Euclidean space. 
\end {thm} 

\begin {thm} \label {evolution-thm2} 
Let the hypotheses of Theorem \ref{evolution-thm1} hold, with $n = \dim(M) = 2$, 
and suppose additionally that $\Omega$ is strictly convex. Let $\partial 
\Omega$ evolve with velocity $k_g \eta$, where $k_g$ is the 
geodesic curvature of $\partial \Omega$. Then 
$$\frac{d}{dt} \log (\lambda) \leq -\frac{4\pi}
{\int_{\partial \Omega} k_g^{-1} d\sigma}.$$
Equality can only occur if $\Omega$ is isometric to a round disk in the 
Euclidean plane.  
\end {thm} 

\begin {thm} \label {evolution-thm3} 
Let the hypotheses of Theorem \ref{evolution-thm1} hold, with $n = \dim (M) \geq 3$, 
and suppose additionally that $\Omega$ is strictly mean convex. Let 
$\partial \Omega$ evolve with velocity $H\eta$, where $H$ is the 
mean curvature of $\partial \Omega$. Then 
$$\frac{d}{dt} \left ( \lambda^{\frac{n-2}{2}} \right )  \leq  
-\left ( \frac{n-2}{2} \right ) \frac{K}
{\int_{\partial \Omega} H^{-1} d\sigma},$$
where $K$ is a constant depending only on $n = \dim (M)$. 
Equality can only occur if $(M,g)$ is 
Euclidean space and $\Omega$ is a round ball. 
\end {thm} 

All three of these theorems are special cases of 
Theorem \ref{evolution-thm4} below. Additionally, one can state similar 
results for other geometric flows, such as flow by Gauss curvature or 
the Willmore flow, under appropriate geometric hypotheses, as corollaries of
Theorem \ref{evolution-thm4}. 

Our secondary goal is a general Schwarz Lemma for the first 
eigenvalue. To make sense of this statement, one should think of 
the Schwarz Lemma geometrically, as interpreted by Pick and 
Ahlfors: if $\mathbf {D}$ is the unit disk in the complex plane 
$\mathbf{C}$, then any holomorphic mapping $F: \mathbf{D} 
\rightarrow \mathbf {D}$ is a contraction in the Poincar\'e metric. 
More recently, Burckel, Marshall, Minda, Poggi-Corradini, and 
Ransford \cite{BMMPR} proves a variety of results in the 
spirit of the Schwarz Lemma for other geometric quantities, such as 
diameter, logarithmic capacity, and area. Bringing the eigenvalue into 
the picture, Laugesen and Morpurgo \cite{LM} prove the following 
result as a special case of a more general theorem: if 
$F:\mathbf{D} \rightarrow \mathbf{C}$ be conformal, then the 
function 
$$r \mapsto \frac{\lambda (F(r\mathbf{D}))}{\lambda (r\mathbf{D})}$$
is a strictly decreasing function, unless $F$ is linear (in which case 
it is constant). Geometrically, the Laugesen-Morpurgo result states that 
$\lambda(F(r\mathbf{D}))$ decreases more rapidly than $\lambda (r 
\mathbf{D})$, as $r$ increases. 

At first glance, all these variations on the Schwarz Lemma seem 
unique to the complex plane. 
However, this is not the case, and one can find several versions of 
the classical Schwarz Lemma in higher dimensions. For instance, 
Yau \cite{Yau} proved a version of the Schwarz Lemma for holomorphic 
mappings between Hermitian manifolds, and later Chen, Cheng, and 
Look \cite{CCL} proved a different version of Yau's result. Motivated 
by these papers, we prove a version of the Schwarz Lemma for the 
first eigenvalue under a conformal diffeomorphism in 
Theorem \ref{schwarz-lemma1} below. 

The rest of the paper is organized as follows. In Section \ref{geometric-sec}
we collect some useful facts from Riemannian geometry. In 
Section \ref{rearrange-sec} we prove some rearrangement and 
reverse-H\"older theorems, which may be of independent interest, 
and finally we prove Theorems \ref{evolution-thm4} and \ref{schwarz-lemma1} 
in Section \ref{main-thm-sec}. 

\bigskip \noindent {\sc Acknowledgements:} We first learned of 
higher-dimensional versions of the Schwarz Lemma from 
Daniel J. F. Fox. We'd like to thank him for the references \cite{Yau} 
and \cite{CCL}. Much of this research was completed 
while J.~R.\  visited University College Cork. He thanks 
UCC for its hospitality. J.~R.\ is partially 
supported by the National Research Foundation of South Africa. 

\section {Geometric preliminaries} \label{geometric-sec}  

In this section we collect some useful results from Riemannian 
geometry. 

\subsection{Model spaces} 

One can find most of the material below in textbooks 
such as \cite{Chav}, or in the survey article \cite{Kar}. 

We begin with some basic formulas regarding our model space $(M_\kappa, 
g_\kappa)$, the complete, simply connected space with constant sectional 
curvature $-\kappa^2$. 
First observe that, by the Cartan-Hadamard theorem, $M$ is diffeomorphic 
to $\R^n$, so we use global polar coordinates. In these 
coordinates, the metric has the form 
\begin {equation} \label{model-metric} 
g_\kappa = dr^2 + \frac{1}{\kappa^2} \sinh^2 (\kappa r) d\theta^2,
\end {equation} where $d\theta^2$ is the round metric on 
the unit sphere. 
%One should be a little careful in interpreting this formula 
%and the one below, as $r$ denotes distance from $0$ in the hyperbolic metric, 
%rather than the Euclidean distance in the more familiar formula for the 
%Poincar\'e metric the unit ball: $g_P = 4(1-\rho^2)^{-2} |dx|^2$. 
Using the expansion 
$$\frac{\sinh(\kappa r)}{\kappa} = r + \frac{1}{3!} \kappa^2 r^3 + \frac{1}{5!}
\kappa^4 r^5 + \cdots,$$
we (formally) recover the Euclidean metric in polar coordinates as $\kappa \rightarrow 
0^+$: 
$$g_0 = dr^2 + r^2 d\theta^2 = \lim_{\kappa \rightarrow 0} \left [ 
dr^2 + \frac{1}{\kappa^2} \sinh^2 (\kappa r)d\theta^2 \right ]. $$
It is also convenient to observe that $\kappa^{-1} \sinh (\kappa r)>
r$ for all $\kappa>0$; geometrically, this says geodesics spread apart 
more rapidly (in fact, exponentially more rapidly) in hyperbolic space 
than in Euclidean space. From \eqref{model-metric} we see that
\begin {equation} \label {model-volumes1}
|\partial B_r|_\kappa = n\omega_n \kappa^{1-n} (\sinh (\kappa r)
)^{n-1} \end {equation} 
and 
\begin {equation} \label{model-volumes2} 
|B_r|_\kappa = n\omega_n \kappa^{1-n} 
\int_0^r (\sinh (\kappa t) )^{n-1} dt = v_\kappa(r),\end {equation}
where $B_r$ is a geodesic ball of radius $r$, and $\omega_n$ 
is the volume of an $n$-dimensional Euclidean unit ball. Later, 
it will be convenient to invert the model volume function $v_\kappa(r)$, 
and write its inverse as $r_\kappa(v)$, which we call the volume radius. 
Again, we can recover the familiar Euclidean formulae by taking a 
limit as $\kappa \rightarrow 0$: 
$$|\partial B_r|_0 = n \omega_n r^{n-1}, \quad |B_r|_0 = \omega_n r^n 
= v_0(r), \quad r_0(v) = \left ( \frac{v}{\omega_n} \right )^{1/n}.$$

The first eigenfunction $\psi_\kappa$ of a geodesic ball in the 
model space $(M_\kappa, g_\kappa)$ is radial, and 
so it satisfies 
\begin {equation} \label{ball-eigen1} 
-\lambda \psi_\kappa = \Delta_\kappa \psi_\kappa = (\sinh (\kappa r))^{1-n} 
((\sinh (\kappa r))^{n-1} \psi_\kappa' )', \end {equation}
where we use $'$ to denote differentiation with respect to $r$. (Where 
it can be understood from context, we suppress the subscript $\kappa$.)
If we change variables to volume and write $\psi^* (v) = \psi(r_\kappa(v))$ 
this equation becomes 
$$-\lambda \psi^* (v) = n^2 \omega_n^2 \kappa^{2-2n} 
\frac{d}{dv} \left ( \sinh^{2n-2}(\kappa r_\kappa(v)) \frac{d\psi^*}{dv} \right ),$$
which we can integrate once to obtain 
\begin {equation} \label{ball-eigen2} 
-(\psi^*)' (v) = n^{-2} \omega_n^{-2} \lambda \left [ \frac{\sinh(\kappa 
r_\kappa (v))}{\kappa} \right ]^{2-2n} \int_0^v \psi^* (t) dt. \end {equation} 

As before, we take a limit as $\kappa \rightarrow 0$ to recover 
the Euclidean analogs of \eqref{ball-eigen1} and 
\eqref{ball-eigen2}, which are (respectively) 
$$-\lambda \psi_0 = \Delta_0 \psi_0 = r^{1-n} (r^{n-1} \psi_0')'$$ 
and 
\begin {equation} \label {ball-eigen3}
-(\psi_0^*)'(v) = n^{-2} \omega_n^{-2/n} \lambda v^{-2 + 2/n} 
\int_0^v \psi_0^* (t) dt.\end {equation}

\subsection{Isoperimetric inequalities} \label{isop-sec} 

In this section we recall some isoperimetric inequalities for general 
Riemannian manifolds. Throughout this section, we take $(M,g)$ 
to be a complete Riemannian manifold, and we usually place a 
bound on its curvature. We also let $\Omega \subset 
M$ be a domain with $\partial \Omega \in \mathcal{C}^\infty$ (though 
this much regularity is rarely necessary), and with $\bar \Omega$ 
compact. A theorem of Beckenbach and Rad\'o \cite{BR} states 
that if $(M,g)$ is a complete surface with nonpositive Gauss curvature 
then 
$$|\partial \Omega|_g^2 \geq 4\pi |\Omega|_g,$$
and equality can only occur if $(M,g)$ is the Euclidean plane and 
$\Omega$ is a round disk. 

The next major break-through is a theorem is due to 
Croke \cite{Croke}, which states that if $\Sect(g) \leq 0$ then 
$$|\partial \Omega|_g \geq c_1(n) |\Omega|_g^{\frac{n-1}{n}},$$
where $c_1(n)$ is an explicit constant and $\Sect(g)$ is the 
sectional curvature of $(M,g)$. This inequality is 
only an equality when $n=4$ and $\Omega$ is a round ball 
in Euclidean space. The next result we quote is due to Kleiner \cite{Klein}, 
and states that if $\dim(M) = 3$ and $\Sect(g) \leq -\kappa^2 \leq 0$ 
then $|\partial \Omega|_g \geq |\partial B_r|_\kappa$, 
where $B_r$ is a geodesic ball in the model space $(M_\kappa, 
g_\kappa)$, with $|\Omega|_g = |B_r|_\kappa$. One only has equality 
in Kleiner's result if $\Omega$ is a model geodesic ball. It is worth remarking
 that Kleiner's proof relies on the Gauss-Bonnet formula, so it 
 can only work in dimension three, while Croke's proof can only be 
 sharp in dimension four. To date, these are the only general results one 
 can find with no restriction on the size of $\Omega$. 
 
 More recently, Morgan and Johnson \cite{MJ} proved a result for compact 
 manifolds $(M,g)$, so long as $|\Omega|_g$ is sufficiently small. Their results 
 state that if $\Sect (g) \leq -\kappa^2$ and $|\Omega|_g$ is sufficiently 
 small then the same inequality $|\partial \Omega|_g \geq |\partial 
 B_r|_\kappa$ holds, where again $B_r$ is a geodesic ball in the 
 model space $(M_\kappa, g_\kappa)$ with $|B_r|_\kappa = 
 |\Omega|_g$. Later, Druet \cite{Druet1} strengthened the Morgan-Johnson 
 result to the point that one only needs a bound on the scalar curvature 
 $S_g$ of $g$, of the form $S_g \leq -n(n-1)\kappa^2$. He also 
 shows that these results hold when $(M,g)$ is complete 
 and $\Omega$ is contained in a small geodesic ball (whose radius 
 might depend on position). Again, one can 
 only have equality in either of these theorems if $\Omega$ is a 
 geodesic ball in the model space. 
 We can use our notation from the previous section to write these inequalities 
 as 
 \begin {equation} \label {isop-small-vol} 
 \Omega \textrm{ sufficiently small } \Rightarrow 
 |\partial \Omega|_g \geq n \omega_n \kappa^{1-n} (\sinh(\kappa 
 r_\kappa (v)))^{n-1}, \end {equation} 
 where $r_\kappa(v)$ is the volume radius function, which 
 inverts \eqref{model-volumes2}. Again, we can contrast this with the 
 Euclidean case, which states that $|\partial \Omega|_0 \geq n \omega_n^{1/n} 
 |\Omega|^{\frac{n-1}{n}}= n\omega_n (r_0(v))^{n-1}$. 
 
 Using these later forms of the isoperimetric inequality for small 
 volumes, Druet \cite{Druet2} and Fall \cite{Fall} proved a 
 Faber-Krahn theorem, and in fact obtained stability estimates. More 
 precisely, they showed that if $\Sect(g) \leq -\kappa^2$ and $|\Omega|_g$ 
 is small, then $\lambda(\Omega) \geq \lambda(B_r)$, 
 where $B_r$ is the geodesic ball in the model space as 
before. Moreover, they estimate the difference $\lambda(\Omega) 
- \lambda(B_r)$, again when $|\Omega|_g$ is small. 

\section {Rearrangements and reverse-H\"older inequalities} 
\label{rearrange-sec}

In this section we discuss a rearrangement of the first eigenfunction
$\phi$ of $\Omega$, and use it to prove an integro-differential inequality 
similar to that of Talenti \cite{Tal}. Next, we obtain  inequalities which 
generalize the results of Chiti \cite{Chiti1, Chiti2} to the Riemannian 
setting. As outlined in our introduction, our standing hypotheses will 
be that \eqref{isop-small-vol} holds. While the precise statements below 
have not yet appeared in the 
literature (to our knowledge), we suspect that much of this section is, 
in the words of A. Treibergs, ``well-known to those who know it well."

Recall that the first eigenfunction $\phi$ satisfies 
\begin {equation} \label {eigenfunction1} 
\lambda(\Omega) = \frac{\int_\Omega |\nabla \phi|^2 \,dm}
{\int_\Omega \phi^2 \,dm} = \inf \left \{ \frac{\int_\Omega 
|\nabla u|^2 \,dm}{\int_\Omega u^2 \,dm} : u \in W^{1,2}_0(\Omega) 
\right \}, \end {equation}
or, alternatively, 
\begin {equation} \label {eigenfunction2} 
\Delta_g \phi  + \lambda(\Omega) \phi = 0, \quad \left. 
\phi\right |_{\partial \Omega} = 0, \quad \phi> 0 
\textrm{ inside }\Omega. \end {equation} 
Let $m = \sup_\Omega \phi$, and for $0 \leq t \leq m$ define 
\begin {equation} \label{dist-funct} 
D_t = \{ \phi > t\}, \qquad \mu(t) = |D_t|_g. 
\end {equation} 
By the co-area formula, we have 
\begin {equation} \label {coarea} 
\mu(t) = \int_t^m \int_{\partial D_\tau} \frac{d\sigma}{|\nabla \phi|}
\,d\tau, \mbox{ so that } \mu'(t) = -\int_{\partial D_t} 
\frac{d\sigma}{|\nabla \phi|} < 0.\end {equation} 
Here $d\sigma$ is the $(n-1)$-dimensional volume element induced on 
$\partial D_t$ by its inclusion in $\Omega$.
Therefore, $\mu$ is monotone, and so it has an inverse 
function we call $\phi^*(v)$, defined by 
$$\phi^*(v) = \inf \{t \in [0,m]: \mu(t) < v\}.$$

While the following is an easy adaptation of equation (34) of \cite{Tal}, 
we include its proof for the reader's convenience. 
\begin {lemma} Let $\Omega \subset (M, g)$ be a domain  with compact closure, 
smooth boundary, and sufficiently small that \eqref{isop-small-vol} holds. 
Then the function $\phi^*$ satisfies 
\begin {equation} \label {eigenfunction3} 
-(\phi^*)'(v) \leq n^{-2} \omega_n^{-2} \lambda(\Omega) 
\left [ \frac{\sinh(\kappa r_\kappa(v))}{\kappa} \right ]^{2-2n} 
\int_0^v \phi^*(t) dt. 
\end {equation} 
Moreover, equality can only occur if $\Omega$ is isometric to a 
geodesic ball in the model space $(M_\kappa, g_\kappa)$. 
\end {lemma} 

\begin {proof} Let $\lambda = \lambda(\Omega)$.  
First observe 
$$\lambda \int_{D_t} \phi dm = 
- \int_{D_t} \Delta_g \phi dm = -\int_{\partial D_t} \frac{\partial \phi}
{\partial \eta} d\sigma = \int_{\partial D_t} |\nabla \phi| d\sigma,$$
 where we have used the divergence theorem and the fact that $\phi$ is 
 constant on $\partial D_t$. Next we use Cauchy-Schwarz to see that 
 $$|\partial D_t|_g \leq \left [ \int_{\partial D_t} |\nabla \phi| d\sigma
\int_{\partial D_t} \frac{d\sigma}{|\nabla \phi|} \right ]^{1/2} = 
\left [ -\lambda \mu'(t) \int_{D_t} \phi dm \right ]^{1/2}.$$ 
Squaring the above inequality and using \eqref{isop-small-vol} will yield 
\begin {equation} \label{integro-diff1} 
-\lambda \mu'(t) \int_{D_t} \phi dm \geq |\partial D_t|_g^2 
\geq |\partial B_r|_\kappa^2 = n^2 \omega_n^2 \kappa^{2-2n} 
(\sinh(\kappa r_\kappa (v)))^{2n-2}. \end {equation} 
Finally we change variables to $v = \mu(t)$, and use the fact that 
$\mu'(t) = \frac{1}{(\phi^*)'(v)}$. This transforms \eqref{integro-diff1} 
into 
$$-(\phi^*)'(v) \leq n^{-2} \omega_n^{-2} \kappa^{2n-2} \lambda
(\sinh(\kappa r_\kappa(v)))^{2-2n} \int_0^v \phi^*(t) dt,$$
 as claimed. \end {proof} 
\noindent  It is vital that \eqref{eigenfunction3} and \eqref{ball-eigen2} have 
 essentially the same right hand side. 
 
\noindent  We also record Talenti's original equality for the eigenfunction, which is 
 \begin {equation} \label {eigenfunction4}
 -(\phi^*)' (v) \leq n^{-2} \omega_n^{-2/n} v^{-2 + 2/n} \int_0^v \phi^*(t) dt. 
 \end {equation}
 This inequality holds when $\Omega \subset (M,g)$ is a domain with compact closure and 
 smooth boundary, with sufficiently small volume, and $S_g \leq 0$. 
 To see \eqref{eigenfunction4}, one can either evaluate a limit 
 of \eqref{eigenfunction3} as $\kappa \rightarrow 0$ or use the same 
 proof with the classical form of the isoperimetric inequality.
 
 We can adapt the arguments of \cite{Chiti1, Chiti2} to our Riemannian 
 setting. 
 \begin {thm} \label{chiti-thm1}
 Let $(M,g)$ and $\Omega$ be as above, so that \eqref{isop-small-vol} 
 applies. Let $B^* \subset (M_\kappa, g_\kappa)$ be a ball in the model space with the same fundamental frequency as $\Omega$, so that $\lambda(B^*) = 
 \lambda(\Omega) = \lambda$.  Let 
 $\phi$ and $\psi$ the the first Dirichlet eigenfunctions of $\Omega$ and of $B^*$ 
 respectively, normalized so that $\|\phi\|_{L^\infty(\Omega)}
 = \|\psi\|_{L^\infty(B^*)}$. Then $\phi^*(v) \geq \psi^*(v)$ for 
 $0 \leq v \leq |B^*|_\kappa$, and equality can only occur for $v>0$ if 
 $\Omega$ is isometric to $B^*$. \end {thm} 
 
 \begin {proof} First observe that, by the Faber-Krahn inequality, we have 
 $|\Omega|_g \geq |B^*|_\kappa$, so that both functions $\psi^*$ and $\phi^*$ 
 are well-defined on the interval $[0,|B^*|_\kappa]$. Moreover, if $|\Omega|_g 
 = |B^*|_\kappa$, then $\Omega$ must be isometric to $B^*$, and we have 
 nothing to prove. Therefore, we can take $|\Omega|_g > |B^*|_\kappa$ 
 without loss of generality. By our normalization of $\psi$, we also know 
 $$\phi^*(0) = \psi^*(0) = m, \quad \psi^*(|B^*|_\kappa) = 0, \quad 
 \phi^* > 0 \textrm{ on }[0,|B^*|_\kappa].$$
 Therefore, there exists $k>1$ such that $k \phi^*(v) > \psi^*(v)$ 
 for all $v \in [0,|B^*|_\kappa]$. Define
 $$k_0 = \inf \{ k > 1 : k \phi^*(v) > \psi^*(v) 
 \textrm{ on }[0,|B^*|_\kappa \}.$$
If $k_0 = 1$ then we've completed the proof, and otherwise there 
exists $v_0 \in (0,|B^*|_\kappa)$ such that $k_0 \phi^*(v_0) 
= \psi^*(v_0)$. If we let 
$$u^*(v) = \left \{ \begin {array}{rl} k_0\phi^*(v) & 0 \leq v \leq v_0 \\ 
\psi^*(v) & v_0 < v \leq |B^*|_\kappa, \end {array} \right. $$
then, by \eqref{eigenfunction3} and \eqref{ball-eigen2}, we have 
\begin {equation} \label{chiti1}
-(u^*)'(v) \leq n^{-2} \omega_n^{-2/n} \lambda \left [ 
\frac{\sinh(\kappa r_\kappa(v))}{\kappa} \right ]^{2-2n} 
\int_0^v u^*(t) dt .\end {equation} 

Now define a radial test function on $B^*$ by $u(r) = u^*(v_\kappa(r))$. We 
use the chain rule and 
$$\frac{dv_\kappa}{dr} = n \omega_n \left ( \frac{\sinh 
(\kappa r)}{\kappa} \right )^{n-1}$$ 
to see that 
\begin {eqnarray*} 
\int_{B^*} |\nabla u|^2 dm & = & \int_0^{|B^*|_\kappa}
n^2 \omega_n^2 \left [ \frac{\sinh (\kappa r_\kappa(v))}
{\kappa} \right ]^{2n-2} (-(u^*)'(v))^2 dv \\ 
& \leq & \lambda \int_0^{|B^*|_\kappa} (-(u^*)'(v)) 
\int_0^v u^*(\tau) d\tau \\ 
& = & \lambda \int_0^{|B^*|_\kappa} u^*(\tau) \int_\tau^{|B^*|_\kappa}
(-(u^*)'(v)) dv d\tau =  \lambda \int_0^{|B^*|_\kappa} (u^*(\tau))^2d\tau \\ 
& = & \lambda \int_{B^*} u^2 dm.\end {eqnarray*}
However, this is impossible unless $u = \psi$, which would contradict 
$k_0 > 1$. 
\end {proof} 

We can integrate the inequality in Theorem \ref{chiti-thm1}
to obtain the following (scale-invariant) corollary. 
\begin {cor} \label{chiti-thm2} 
Let $\Omega \subset (M,g)$ be as above, and let $B^*$ 
be the geodesic ball in the model space $(M_\kappa, g_\kappa)$ 
with $\lambda(\Omega) = \lambda(B^*)$. Let $\psi$ be the 
first eigenfunction of $B^*$ and let $\phi$ be the first eigenfunction 
of $\Omega$. Then for all $p>0$ we have 
$$\frac{\|\phi\|_{L^p(\Omega)}}{\|\phi\|_{L^\infty(\Omega)}} 
\geq \frac{\|\psi\|_{L^p(B^*)}}{\|\psi\|_{L^\infty(B^*)}}.$$
Equality can only occur if $\Omega$ is isometric to $B^*$. 
\end {cor} 

One can find a version of the following theorem, which reverses the 
standard H\"older inequality,  in the hyperbolic setting  
in Section 9 of \cite{BL}. Both proofs utilize Chiti's method 
from \cite{Chiti2}. 

\begin {thm} With the same $\Omega$ as above and any choice 
$0 < p < q < \infty$, there exists a 
positive, finite constant $C = C(n,p,q,\kappa, \lambda)$ such that 
the first eigenfunction $\phi$ of $\Omega$, with eigenvalue $\lambda$, 
satisfies 
\begin {equation} \label{reverse-holder1} 
\left ( \int_\Omega \phi^p dm \right )^q \geq C \left ( 
\int_\Omega \phi^q dm \right )^p. \end {equation} 
Equality can only occur if $\Omega$ is isometric to 
$B^*$. 
\end {thm} 
In fact, it will be transparent from the proof that 
\begin {equation} \label{holder-const1}
C = \frac{\left ( \int_{B^*} \psi^p dm \right )^q}
{\left ( \int_{B^*} \psi^q dm \right )^p},\end {equation}
where $B^*$ is the geodesic ball in the model space 
with $\lambda(B^*) = \lambda = \lambda(\Omega)$, 
and $\psi$ is its first eigenfunction. 

\begin {proof} We use the same approach as in 
the proof of Theorem \ref{chiti-thm1}, but this time normalize $\psi$
such that 
$$\int_{B^*}\psi^p dm = \int_\Omega \phi^p dm.$$
Thus, by Corollary \ref{chiti-thm2} above, $\| \psi\|_{L^\infty(B^*)} \geq 
\| \phi\|_{L^\infty(\Omega)}$, with equality if and only if $\Omega$ 
is isometric to $B^*$. We may therefore assume 
\begin {equation} \label{rearrange1} 
\psi^*(0) = \| \psi\|_{L^\infty(B^*)} >\|\phi\|_{L^\infty(\Omega)} = \phi^*(0).
\end {equation} 
We also know, as before, that 
$$\psi^*(|B^*|_\kappa) = 0, \quad \phi^* > 0 \textrm{ on } 
[0, |B^*|_\kappa],$$
which combined with \eqref{rearrange1} tells us the graphs 
of $\phi^*$ and $\psi^*$ must cross, and not just touch, at least once on the interval 
$[0,|B^*|_\kappa]$. Define 
$$v_0 = \sup \{ v \in (0, |B^*|_\kappa) : \phi^*(\tilde v) \leq \psi^*(\tilde v) 
\textrm{ for all }\tilde v \in (0,v) \}, $$
so that we have 
$$0 < v_0 < |B^*|_\kappa, \quad \psi^* \geq \phi^* \textrm{ in }[0, v_0], 
\quad \phi^*(v_0) = \psi^*(v_0). $$
Additionally, there must exist $\delta>0$ such that $\phi^*(v) > \psi^*(v)$ 
for $v \in (v_0, v_0 + \delta)$. 

We claim that actually $\phi^* > \psi^*$ in the interval $(v_0, |B^*|_\kappa]$. 
Indeed, if this were not the case then there would exist $v_1$
such that 
$$v_0 < v_1 < |B^*|_\kappa, \quad \psi^*(v_1) = \phi^*(v_1), 
\quad  \phi^*(v) > \psi^*(v) \textrm { for }v_0 < v< v_1.$$ 
This allows us to define a test function for $B^*$ as 
$$u^*(v) = \left \{ \begin {array} {rl} \psi^*(v) & 0 \leq v \leq v_0 \\ 
\phi^*(v) & v_0 \leq v \leq v_1 \\ \psi^* (v) & v_1 \leq v \leq |B^*|_\kappa . 
\end {array} \right . $$
As before, our test function satisfies \eqref{chiti1}, 
$$-(u^*)'(v) \leq n^{-2} \omega_n^{-2/n} \lambda \left [ 
\frac{\sinh(\kappa r_\kappa(v))}{\kappa} \right ]^{2-2n} 
\int_0^v u^*(t) dt,$$
and we can define a radial test function $u$ on $B^*$ by 
$u(r) = u^* (v_\kappa(r))$, which in turn satisfies 
\begin {eqnarray*} 
\int_{B^*} |\nabla u|^2 dm & = & \int_0^{|B^*|_\kappa}
n^2 \omega_n^2 \left [ \frac{\sinh (\kappa r_\kappa(v))}
{\kappa} \right ]^{2n-2} (-(u^*)'(v))^2 dv \\ 
& \leq & \lambda \int_0^{|B^*|_\kappa} (-(u^*)'(v)) 
\int_0^v u^*(\tau) d\tau \\ 
& = & \lambda \int_0^{|B^*|_\kappa} u^*(\tau) \int_\tau^{|B^*|_\kappa}
(-(u^*)'(v)) dv d\tau =  \lambda \int_0^{|B^*|_\kappa} (u^*(\tau))^2d\tau \\ 
& = & \lambda \int_{B^*} u^2 dm.\end {eqnarray*}
As before, this is only possible if $u= \psi$, which contradicts 
our assumption $\psi^*(0) > \phi^*(0)$. 

So far, we have shown there exists $v_0 \in (0, |B^*|_\kappa)$ 
such that $\psi^*\geq\phi^*$ on $(0,v_0)$ and $\phi^* > \psi^*$ 
on $(v_0, |B^*|_\kappa)$. We extend $\psi^*$ to be zero on the 
interval $[|B^*|_\kappa, |\Omega|_g]$, and claim that 
\begin {equation} \label{chiti2} 
v\in [0, |\Omega|_g] \Rightarrow 
\int_0^v (\psi^*(\tau))^p d\tau \geq \int_0^v (\phi^*(\tau))^p d\tau. 
\end {equation} 
To prove this claim, we let 
$$I(v) = \int_0^v (\psi^*(\tau))^p d\tau - \int_0^v 
(\phi^*(\tau))^p d\tau$$ 
and observe 
$$I(0) = I(|\Omega|_g) = 0, \quad  I' (v) = (\psi^*(v))^p 
- (\phi^*(v))^p.$$
Thus $I$ is increasing on the interval $[0, v_0)$ and decreasing 
on the interval $(v_0, |\Omega|_g]$. It follows immediately that 
$I(v) > 0$ for $0 \leq v \leq v_0$. If we had $I(v_1) < 0$ for some 
$v_1 \in (v_0, |\Omega|_g)$ then, because $I$ is decreasing 
in this interval, we would also have $I(|\Omega|_g) < 0$, which is a 
contradiction. We conclude \eqref{chiti2}. It follows from an 
inequality of Hardy, Littlewood, and P\'olya \cite{HLP} that for all $q>p$ 
we have 
$$\left ( \int_\Omega \phi^q dm \right )^{1/q} \leq \left ( \int_{B^*} \psi^q 
dm \right )^{1/q} = \frac{\left ( \int_{B^*} \psi^q dm \right )^{1/q} }
{\left ( \int_{B^*} \psi^p dm \right )^{1/p}} \cdot \left ( \int_{B^*}
\phi^p dm \right )^{1/p},$$
which we can rearrange to read 
$$\frac{\left (\int_\Omega \phi^p dm \right )^{1/p}}
{\left ( \int_\Omega \phi^q dm \right )^{1/q}} \geq \frac
{\left ( \int_{B^*} \psi^p dm \right )^{1/p}}{\left (\int_{B^*}
\psi^q dm \right )^{1/q}} = \tilde C.$$
Raising this inequality to the power $pq$, we then obtain
\[
\left ( \int_\Omega \phi^p dm \right )^q \geq \tilde C^{pq} 
\left (\int_\Omega \phi^q dm  \right )^p.
\qedhere\]
\end {proof} 

In the case $S_g \leq 0$ we can extract the explicit dependence of the constant 
$C$ in \eqref{reverse-holder1} on the eigenvalue $\lambda$. The 
dependence on the eigenvalue in the hyperbolic case is more 
challenging to understand, because the eigenfunctions 
on geodesic balls do not scale in curved setting (see, for 
instance, Section 3 of \cite{BL}). 
\begin {cor} Suppose $\Omega$ is a  domain in $(M,g)$, where $S_g \leq 0$, 
which is sufficiently small so that \eqref{isop-small-vol} applies. 
Let $\phi$ be its first  eigenfunction, with eigenvalue $\lambda$. 
Then there is a constant $K = K(n,p,q)$ such that 
\begin {equation} \label{reverse-holder2} 
\left ( \int_\Omega \phi^p dm \right )^q \geq K 
\lambda^{n(p-q)/2} \left ( \int_\Omega \phi^q dm \right )^p. 
\end {equation} 
\end {cor} 

\begin {proof} This time our comparison domains are round 
balls in Euclidean space, and the dilation of an eigenfunction 
on a ball is an eigenfunction on the corresponding dilated ball. 
We have, according to \eqref{holder-const1}
$$C = \frac{\left ( \int_{B^*} \psi^p dm \right )^q}
{\left ( \int_{B^*} \psi^q dm \right )^p}.$$
Denote the Euclidean  radius of $B^*$ by $\rho$, and change variables 
to the unit ball by defining the function $\tilde \psi(r) 
= \psi (r\rho)$, so that 
$$C = \rho^{n(q-p)} \frac{\left (\int_{B_1} \tilde \psi^p \,dm 
\right )^q} {\left (\int_{B_1} \tilde \psi^q \,dm \right )^p}.$$
Now, $\tilde \psi$ is the first eigenfunction on the unit 
ball in Euclidean space, and all that remains is to recall the 
scaling law for eigenvalues: $\lambda (B^*) = \rho^{-2} \lambda(B_1)$. 
Thus we see that
$$\rho = \left ( \frac{\lambda (B^*)}{\lambda(B_1)} 
\right )^{-1/2} = \left ( \frac{\lambda (\Omega)}{\lambda(B_1)}
\right )^{-1/2},$$
and so 
\[
C = \lambda^{-\frac{n}{2} (q-p)} \lambda(B_1)^{\frac{n}{2} 
(q-p)} \frac{\left ( \int_{B_1} \tilde \psi^p \,dm \right )^q}
{\left ( \int_{B_1} \tilde \psi^q \,dm \right )^p}.
\qedhere
\]
\end {proof} 
We will later use the case of $p=1$ and $q=2$, which reads 
\begin {equation} \label{reverse-holder3} 
\left ( \int_\Omega \phi \,dm \right )^2 \geq K \lambda^{-n/2} 
\int_\Omega \phi^2 \,dm.\end {equation} 
In the case of $\dim(M) = 2$ we recover an inequality of 
Payne and Rayner \cite{PR}: 
\begin {equation} \label{reverse-holder4} 
\left ( \int_\Omega \phi \,dm \right )^2 \geq \frac{4\pi}{\lambda}
\int_\Omega \phi^2 \,dm.
\end {equation} 
Here we have used the sharp version of the isoperimetric inequality 
of Beckenbach and Rad\'o \cite{BR} for complete surfaces with 
nonpositive Gauss curvature. It is also important to notice that in the 
two-dimensional case we do not place any restriction on the size 
of $\Omega$. 

The reverse Cauchy-Schwarz inequality \eqref{reverse-holder4} can be rewritten as a geometric 
isoperimetric inequality for the (singular) conformal metric 
$\tilde g = |\nabla \phi|^2 g$. We have the following corollary. 
\begin {cor} Let $(M,g)$ is a surface with nonpositive Gauss curvature, 
and let $\Omega$ be a domain with $\bar \Omega$ compact and 
$\partial \Omega \in \mathcal{C}^\infty$. Place the (singular) conformal
metric $\tilde g = |\nabla \phi|^2 g$ on $\Omega$, where $\phi$ is the 
first Dirichlet eigenfunction of $\Delta_g$ on $\Omega$. Then, with 
respect to $\tilde g$, we have 
$$\tilde L^2 \geq 4\pi \tilde A,$$
and equality can only occur if $\Omega$ is isometric to a flat 
disk. 
\end {cor}

\begin {proof} We begin with the left hand side of \eqref{reverse-holder4}. We 
have 
$$\left ( \int_\Omega \phi dm \right )^2 = \frac{1}{\lambda^2} \left ( \int_\Omega 
\Delta \phi dm \right )^2 = \frac{1}{\lambda^2} \left ( \int_{\partial \Omega} 
\frac{\partial \phi}{\partial \eta} d\sigma \right )^2 = \frac{1}{\lambda^2} 
\left ( \int_{\partial \Omega} |\nabla \phi| d\sigma \right )^2 = \frac{\tilde L^2}
{\lambda^2},$$
where we have used the PDE satisfied by $\phi$, the divergence theorem, and 
the fact that $\phi$ is constant on $\partial \Omega$. 
On the other hand, the right hand side of \eqref{reverse-holder4} is 
$$\frac{4\pi}{\lambda} \int_\Omega \phi^2 dm = \frac{4\pi}{\lambda^2} 
\int_{\Omega} |\nabla \phi|^2 dm = \frac{4\pi \tilde A}{\lambda^2}.$$
The result follows. \end {proof} 

\section {Monotonicity of the first eigenvalue} \label{main-thm-sec}

In this section we study the evolution of $\lambda$ as 
$\Omega$ evolves.

A key ingredient is the reverse-H\"older inequality for the first eigenfunction
we developed in Section \ref{rearrange-sec}. Another key 
ingredient is, naturally, the Hadamard variation fomula for the first eigenvalue. 
We consider a one-parameter family of 
diffeomorphisms $\zeta (t,p) : (-\epsilon, \epsilon) \times M \rightarrow 
M$, and let $\Omega_t = \zeta(t,\Omega)$. The family of 
mappings $\zeta$ is the flow of the time-dependent vector field 
$\chi$, where 
\begin {equation} \label{flow-eqn}
\frac{\partial \zeta}{\partial t} (t,p) = \chi (t,p).\end {equation}
In this way, if $\Omega = \Omega_0$ satisfies our standing 
hypotheses, then so will $\Omega_t$ for $t$ sufficiently small. 

We let $\lambda(t) = \lambda(\Omega_t)$, and use a dot to denote 
differentiation with respect to $t$. A classical theorem of 
Hadamard \cite{Had} states that
\begin {equation} \label {had-var} 
\dot \lambda(0) = - \int_{\partial \Omega} \langle \chi, \eta 
\rangle \left ( \frac{\partial \phi}{\partial \eta} \right )^2 d\sigma, 
\end {equation}
where $\phi$ is the first eigenfunction of $\Omega$, normalized so 
that $\int_\Omega \phi^2 dm = 1$. We include the proof 
for the reader's convenience. 
\begin {proof} First we compute the time derivative 
of the boundary terms 
of the normalized first eigenfunction $\phi$. Taking 
a derivative of the condition 
\[
\phi\big(t,\zeta(t,p)\big) = 0,\  p \in 
\partial \Omega
\]
with respect to $t$ and using \eqref{flow-eqn}, we 
obtain
\[
\dot\phi\big(t,\zeta(t,p)\big) + 
\big\langle \nabla\phi\big(t, \zeta(t,p)\big),
\chi(p) \big\rangle = 0.
\]
Here and later, the gradient refers only to the 
spatial derivative. Set $t=0$ and use the fact 
that $\phi$ is constant along $\partial \Omega_t$ to 
obtain 
\begin {equation} \label {first-var-a} 
\dot \phi(0,p)  
= - \big\langle \nabla \phi (0,p), \chi(p) 
\big\rangle = -\Big\langle \left.\frac{\partial \phi} 
{\partial \eta} \right|_{(0,p)} \eta(p), \chi(p)
\Big\rangle, \quad p \in \partial \Omega.
\end {equation} 

Next we take the derivative of the eigenfunction 
equation 
\begin {equation} \label {first-var-b} 
\Delta \phi\big(t,\zeta(t,p)\big) + \lambda(t) \phi\big(t,\zeta(t,p)\big) = 0 
\end {equation} 
with respect to $t$, leading to
\begin {eqnarray*}
0 & = & \Delta \left [ \dot \phi + \langle \nabla 
\phi, \chi \rangle \right] + \lambda(t) 
\left [ \dot \phi + \langle \nabla \phi, 
\chi\rangle \right ] + \dot \lambda(t) \phi\\
& = & \Delta \dot\phi + \langle \nabla \Delta \phi, 
\chi \rangle  +  \lambda(t) 
\dot\phi + \lambda(t) \langle \nabla \phi, 
\chi\rangle + \dot \lambda(t) \phi \\
& = & \Delta \dot\phi + \lambda(t) \phi_t + 
\dot\lambda(t) \phi .
\end {eqnarray*}
Setting $t=0$ and rearranging yields 
\begin {equation} \label{first-var-c} 
\Delta \left. \dot\phi \right|_{t=0} + \lambda(0) 
\left. \dot\phi \right|_{t=0} 
= -\dot\lambda(0)  \phi\big\vert_{t=0} \quad 
\text{ in } \Omega.
\end {equation} 
We multiply  (\ref{first-var-b}), with $t=0$,  
by $\dot\phi\big\vert_{t=0}$ and 
multiply (\ref{first-var-c}) by $\phi$, subtract and 
obtain 
\begin {equation} \label {first-var-d}
\dot\lambda(0) \phi^2(0,p) = \dot\phi(0,p) 
\Delta \phi(0,p) - \phi(0,p) 
\Delta \dot\phi(0,p), \quad p \in \Omega.
\end {equation} 

Integrate \eqref{first-var-d} over $\Omega$ and 
use the fact that $\int_{\Omega_t} \phi^2 dm  = 1$ 
to obtain 
\begin {eqnarray*} 
\dot\lambda(0) & = & \int_{\Omega} \dot\phi \Delta \phi 
- \phi \Delta \dot\phi dm\\
& = & \int_{\partial \Omega} \dot\phi 
\frac{\partial \phi}{\partial \eta} d\sigma- 
\int_{\Omega} \big\langle \nabla \phi, 
\nabla \dot\phi \big \rangle d\sigma+ 
\int_{\Omega} \big\langle \nabla \phi, 
\nabla \dot\phi \big\rangle d\sigma- \int_{\partial \Omega} 
\phi \frac{\partial \dot\phi}{\partial \eta} d\sigma \\
& = &\int_{\partial \Omega} \dot\phi 
\frac{\partial \phi}{\partial \eta} d\sigma \\
& = & - \int_{\partial \Omega} \frac{\partial \phi}
{\partial\eta} \langle \nabla \phi, \chi \rangle d\sigma\\
& = & - \int_{\partial \Omega} \langle \chi, 
\eta\rangle \left ( \frac{\partial \phi}{\partial \eta} 
\right )^2 d\sigma, 
\end {eqnarray*} 
which is equation (\ref{had-var}) as claimed. In the second equality above we 
integrated by parts, in the next to last we used \eqref{first-var-a}, and at the last step 
we used the fact that $\phi$ is constant on $\partial \Omega$ (and hence $\nabla \phi = 
\frac{\partial \phi}{\partial \eta} \eta$ there). 
\end {proof} 

We will need to transform \eqref{reverse-holder3} and \eqref{reverse-holder4} 
for our use later in bounding $\dot \lambda$. 
\begin {lemma}  Let $(M,g)$ be a complete Riemannian manifold with 
nonpositive scalar curvature, and let $\Omega$ be a sufficiently small 
domain in $M$ with $\bar \Omega$ compact and $\partial \Omega \in 
\mathcal {C}^\infty$, so that \eqref{isop-small-vol} applies. Let $\phi$ 
be the first Dirichlet eigenvalue of $\Delta_g$ on $\Omega$, normalized 
so that $\int_\Omega \phi^2 dm = 1$. Then 
\begin {equation} \label {reverse-holder5} 
K \lambda^{2-\frac{n}{2}} \leq \left ( \int_{\partial \Omega}
\frac{\partial \phi}{\partial \eta} d\sigma \right )^2,
\end {equation} 
where $K$ is the same constant, depending only on $n$, 
in \eqref{reverse-holder3}. In dimension two, this inequality reads
 \begin {equation} \label {reverse-holder6} 
 4\pi \lambda \leq \left ( \int_{\partial \Omega} \frac{\partial \phi}
 {\partial \eta} d\sigma \right )^2. \end {equation} 
 Equality can only occur if $\Omega$ is isometric to a flat ball in the 
 appropriate dimensional Euclidean space. 
 \end {lemma} 

\begin {proof} By our normalization we have 
\begin {eqnarray*}
K \lambda^{-n/2} & = & K \lambda^{-n/2} \int_\Omega \phi^2 dm 
\leq \left ( \int_\Omega \phi dm \right )^2 \\
& = & \frac{1}{\lambda^2} 
\left ( \int_\Omega \Delta \phi dm \right )^2 = \frac{1}{\lambda^2}
\left ( \int_{\partial \Omega}  \frac{\partial \phi}{\partial \eta} 
d\sigma \right )^2.\end {eqnarray*}
The result follows. \end {proof} 

Recall that we have set $\Omega_t = \zeta(t,\Omega)$, where 
$\zeta$ is a one-parameter family of diffeomorphisms on $M$. We 
let $\lambda(t) = \lambda(\Omega_t)$, and we are assuming all 
the hypotheses relevant to \eqref{isop-small-vol} hold. 

\begin {thm} \label{evolution-thm4}
Let $\partial \Omega$ move with velocity $e^{w} \eta$, where 
$\eta$ is the unit outward normal of $\partial \Omega$ and $w$ is 
a bounded and continuous function. If $n = \dim (M) \geq 3$ then 
\begin {equation} \label {evolution1} 
\frac{d}{dt} \left [ \lambda^{\frac{n-2}{2}} \right ] 
\leq - \left ( \frac{n-2}{2} \right ) \frac{K}{\int_{\partial \Omega} 
e^{-w} d\sigma}, \end {equation} 
where $K$ is the same constant in \eqref{reverse-holder3}, which 
depends only on $n$. If $\dim (M) =2$ then 
\begin {equation} \label{evolution2} 
\frac{d}{dt} \log (\lambda) \leq -\frac{4\pi}{\int_{\partial \Omega} 
e^{-w} d\sigma}. \end {equation} 
In either case, equality can only occur if $\Omega$ is 
isometric to a round ball in the appropriate dimensional Euclidean 
space. 
\end {thm} 

\begin {proof} By Cauchy-Schwarz, 
$$-\int_{\partial \Omega} \frac{\partial \phi}{\partial \eta} 
d\sigma = \int_{\partial \Omega} \left ( e^{-w/2} \right ) \left ( e^{w/2}
\left | \frac{\partial \phi}{\partial \eta} \right | \right ) d\sigma \leq \left ( 
\int_{\partial \Omega} e^{-w} d\sigma \right )^{1/2} \left ( \int_{\partial \Omega} 
e^w \left ( \frac{\partial \phi}{\partial \eta} \right )^2 d\sigma \right )^{1/2},$$
 so that 
 \begin {equation} \label{cauchy-schwarz1} 
 \int_{\partial \Omega} e^w \left ( \frac{\partial \phi}{\partial \eta} 
 \right )^2 d\sigma \geq \frac{1}{\int_{\partial \Omega} e^{-w} d\sigma} 
 \left ( \int_{\partial \Omega} \frac{\partial \phi}{\partial \eta} d\sigma 
 \right )^2. \end {equation} 
We first prove \eqref{evolution1}. 
Using \eqref{had-var}, \eqref{cauchy-schwarz1},  
and \eqref{reverse-holder5}, we see 
\begin {eqnarray*}
-\dot \lambda & = & \int_{\partial \Omega} e^{w} \left (
\frac{\partial \phi}{\partial \eta} \right )^2 d\sigma \\ 
& \geq & \frac{1}{\int_{\partial \Omega} e^{-w} d\sigma}
\left ( \int_{\partial \Omega}  \frac{\partial \phi}{\partial \eta}
d\sigma \right )^2 \\ 
& \geq & \frac{K \lambda^{2-\frac{n}{2}}}{\int_{\partial \Omega}
e^{-w} d\sigma }, \end {eqnarray*}
which we can rearrange to read 
$$- \frac{d}{dt} \left [ \frac{2}{n-2} \lambda^{\frac{n-2}{2}} 
\right ] = - \lambda^{\frac{n}{2} - 2} \dot \lambda \geq 
\frac{K}{\int_{\partial \Omega} e^{-w} d\sigma}.$$

The proof of \eqref{evolution2} is very similar. This time we 
replace \eqref{reverse-holder5} with \eqref{reverse-holder6} 
to obtain 
\begin {eqnarray*} 
-\dot \lambda & = & \int_{\partial \Omega} e^w \left ( 
\frac{\partial \phi}{\partial \eta} \right )^2 \\ 
& \geq & \frac{1}{\int_{\partial \Omega} e^{-w} d\sigma} 
\left ( \int_{\partial \Omega} \frac{\partial \phi}{\partial \eta}
d\sigma \right )^2 \\ 
& \geq & \frac{4\pi \lambda }{\int_{\partial \Omega} e^{-w} 
d\sigma} ,\end {eqnarray*} 
which we can rearrange to read 
$$- \frac{d}{dt} \log \lambda = -\frac{\dot \lambda}{\lambda} 
\geq \frac{4\pi}{\int_{\partial \Omega} e^{-w} d\sigma}.$$
\end {proof} 

Now Theorem \ref{evolution-thm1} follows by taking $w=0$,
Theorem \ref{evolution-thm2} follows by taking $w= \log{k_g}$, and 
Theorem \ref{evolution-thm3} follows by taking $w= \log{H}$. 

Finally, we apply our technique to the case that $\Omega$ is the 
conformal image of a Euclidean ball. We let $(M,g)$ be a complete 
Riemannian manifold of dimension $n$ with $S_g \leq 0$ and 
let $F:\R^n \rightarrow M$ be a conformal mapping. Let $B_t$ 
be the ball of radius $t$ in $\R^n$, and let $\Omega_t = F(B_t)$. Letting 
$\lambda(t) = \lambda(B_t)$ and $\tilde \lambda(t) = \lambda (F(B_t))$, 
we wish to compare $\lambda(t)$ to $\tilde \lambda(t)$.  As $t$ increases, 
$\partial B_t$ moves with velocity $\eta= \frac{\partial}{\partial r}$, and 
(because $F$ is conformal) $\partial \Omega$ moves with velocity $|DF|
\tilde \eta$, where $\tilde \eta$ is the outward unit normal of $\Omega_t$. 
We let $\phi$ be the first Dirichlet eigenfunction of $\Delta$ on $B_t$, 
normalized so that $\int_{B_t} \phi^2 dm = 1$, and let $\tilde \phi$  be 
the first Dirichlet eigenfunction of $\Delta_g$ on $\Omega_t$, normalized 
so that $\int_{\Omega_t} \tilde \phi^2 dm = 1$. It will be convenient to define 
$\psi = \tilde \phi \circ F$, and observe that $|\nabla \psi| = |DF| |\nabla \tilde \phi|$. 

\begin {thm} \label{schwarz-lemma1}
Let $F:\R^n \rightarrow M$ be conformal, where $(M,g)$ is complete, 
with $S_g \leq 0$ as above. If $n=2$ then 
\begin {equation} \label {schwarz-ineq1}
\frac{d}{dt} \log (\tilde \lambda/\lambda) < 0 \end {equation}
unless $F$ is an isometry when restricted to $B_t$. If $n \geq 3$,
$t$ is small enough so that \eqref{isop-small-vol} applies to $\Omega_t$, 
and $\int_{\partial B_t} |DF|^{n-2} d\sigma > |\partial B_t| = n \omega_n t^{n-1}$ 
then 
\begin {equation} \label{schwarz-ineq2}
\frac{d}{dt} \left [ \tilde \lambda^{\frac{n-2}{2}} - 
\lambda^{\frac{n-2}{2}} \right ] < 0. \end {equation} 
\end {thm} 

Notice that we recover the (a special case of) the Laugesen-Morpurgo result 
in \cite{LM} in dimension two. In higher dimensions, this theorem states that 
if $F$ is a conformal map with a sufficiently large coformal factor then 
$\tilde \lambda^{\frac{n-2}{2}}$ decreases more rapidly than 
$\lambda^{\frac{n-2}{2}}$. Thus, our theorem is very much in the 
spirit of the results in \cite{LM} and in \cite{BMMPR}.

\begin {proof} First observe that, because $\partial \Omega_t$ 
moves with velocity $|DF|\tilde \eta$, the Hadamard variation formula 
becomes 
\begin {equation} \label{had-var2}
-\dot {\tilde \lambda} = \int_{\partial \Omega_t} |DF| \left ( \frac{\partial 
\tilde \phi}{\partial \tilde \eta} \right )^2 d\tilde \sigma = \int_{\partial B_t}
|DF|^{n-2} \left ( \frac{\partial \psi}{\partial \eta} \right )^2d\sigma.
\end {equation} 
Thus, in dimension $n \geq 3$, the inequality \eqref{reverse-holder5} 
gives 
\begin {eqnarray*}
K \tilde \lambda^{\frac{4-n}{2}} & \leq &\left ( \int_{\partial 
\Omega_t} |\nabla \tilde \phi | d\tilde \sigma \right )^2\\ 
& = & \left ( \int_{\partial B_t} |DF|^{n-2} |\nabla \psi| d\sigma 
\right )^2 \\ 
& \leq & \int_{\partial B_t} |DF|^{n-2} d\sigma  \cdot 
\int_{\partial B_t} |DF|^{n-2} |\nabla \psi|^2 d\sigma  \\ 
& = & - \dot{\tilde \lambda} \int_{\partial B_t} |DF|^{n-2} d\sigma, 
\end {eqnarray*} 
which we can rearrange to give 
$$-\frac{d}{dt} \left [ \frac{2}{n-2}  \tilde \lambda^{\frac{n-2}{2}} \right ]
= -\tilde \lambda^{\frac{n-4}{2}} \dot {\tilde \lambda} \geq 
\frac{K}{\int_{\partial B_t} |DF|^{n-2} d\sigma}.$$
However, the equality case in Theorem \ref{evolution-thm1} tells us 
$$-\frac{d}{dt} \left [ \frac{2}{n-2} \lambda^{\frac{n-2}{2}} \right ] 
= \frac{K}{|\partial B_t|},$$
so \eqref{schwarz-ineq2} now follows from the inequality $\int_{\partial 
\Omega} |DF|^{n-2} d\sigma  > |\partial B_t|$. In the two-dimensional 
case, we use \eqref{reverse-holder6} to see 
\begin {eqnarray*} 
4\pi \tilde \lambda & \leq & \left ( \int_{\partial \Omega_t} 
|\nabla \tilde \phi| d\tilde \sigma \right )^2  =  \left ( 
\int_{\partial B_t}  |\nabla \psi| d\sigma \right )^2 \\
& \leq & |\partial B_t| \int_{\partial B_t} |\nabla \psi|^2 d\sigma \\
& = & -\dot {\tilde \lambda} |\partial B_t|,  \end {eqnarray*}
which we can rearrange to give 
$$-\frac{d}{dt} \log \tilde \lambda = -\frac{\dot {\tilde \lambda}}{\tilde \lambda} \geq 
\frac{4\pi}{|\partial B_t|} = - \frac{\dot \lambda}{\lambda} = -\frac{d}{dt}
\log \lambda,$$
where we have again used the equality case of Theorem \ref{evolution-thm1}. This 
completes the proof of \eqref{schwarz-ineq1}. 
\end {proof}

\end{document}